\documentclass{amsart}
\usepackage{amsmath}
\usepackage{amssymb}
\usepackage{latexsym}
\usepackage{verbatim}
\usepackage{graphics}
\usepackage{amsfonts}
\usepackage{graphicx}
\usepackage{epstopdf}
\newtheorem{thm}{Theorem}
\newtheorem{conj}{Conjecture}
\newtheorem{prop}{Proposition}
\newtheorem{lem}{Lemma}
\newtheorem{defn}{Definition}
\newtheorem{cor}{Corollary}
\newtheorem{quest}{Question}

\begin{document}

\title{Counting Links in Complete Graphs}
\date{\today}

\author{Thomas Fleming}
\address{Department of Mathematics\\
   	University of California, San Diego\\
   	La Jolla, CA 92093-0112}
\email{tfleming@math.ucsd.edu}

\author{Blake Mellor}
\address{Mathematics Department\\
   	Loyola Marymount University\\
   	Los Angeles, CA  90045-2659}
\email{bmellor@lmu.edu}


\begin{abstract}

We find the minimal number of links in an embedding of any complete $k$-partite graph on 7 vertices (including $K_7$, which has at least 21 links).  We give either exact values or upper and lower bounds for the minimal number of links for all complete $k$-partite graphs on 8 vertices.  We also look at larger complete bipartite graphs, and state a conjecture relating minimal linking embeddings with minimal book embeddings.

\end{abstract}

\maketitle



\section{Introduction} \label{S:intro}

\begin{figure} [h]
    $$\includegraphics{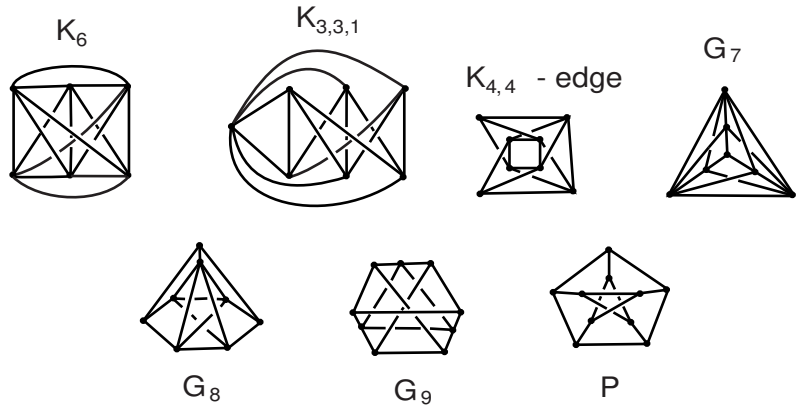}$$
    \caption{The Petersen family of graphs} \label{F:petersen}
    \end{figure}

The study of links and knots in spatial graphs began with Conway and Gordon's seminal result that every embedding of $K_6$ contains a non-trivial link and every embedding of $K_7$ contains a non-trivial knot \cite{cg}.  Their result sparked considerable interest in {\it intrinsically linked} and {\it intrinsically knotted} graphs -- graphs with the property that every embedding in $\mathbb{R}^3$ contains a pair of linked cycles (respectively, a knotted cycle).  Robertson, Seymour and Thomas \cite{rst} gave a Kuratowski-type classification of intrinsically linked graphs, showing that every such graph contains one of the graphs in the {\it Petersen family} as a minor (see Figure~\ref{F:petersen}).

However, while their result answers the question of {\it which} graphs are intrinsically linked, it does not address {\it how} they are linked, and how complicated the linking must be.  In this paper, we measure the ``complexity" of a graph (with respect to intrinsic linking) by the minimal number of links in any embedding of the graph (denoted $mnl(G)$).

This is not the only possible measure of complexity.  Rather than counting the number of links, one could focus instead on the complexity of the individual links.  Flapan \cite{fl} has given examples of graphs which must contain links with large linking numbers, and Flapan et al. \cite{ffnp} constructed graphs whose embeddings must contains links with many components.  Recently, the second author, with Flapan and Naimi, has generalized these results to show that there are graphs whose embedding must contain a link with arbitrarily complex {\it linking patterns} \cite{fmn}.

In other work, the authors have used the notion of {\it virtual spatial graphs} to form a filtration of graphs based on the presence of virtual links in a graph's virtual diagrams.  These various measures of complexity are quite different; for example, while all the graphs in the Petersen family have $mnl(G) = 1$, in the virtual filtration they fall into two distinct levels.  While both $K_6$ and $K_7$ can be embedded with only Hopf links, $mnl(K_6) = 1$ while we will show that $mnl(K_7) = 21$.

Our goal in this paper is to count the minimal number of links in small (7 or 8 vertices) complete $k$-partite graphs.  We obtain complete results for graphs on 7 vertices, and upper and lower bounds for $mnl(G)$ for graphs on 8 vertices (see Table \ref{Ta:results}).  In the final sections, we look at larger complete bipartite graphs, and conjecture a relationship between minimal number of links and minimal book embeddings of graphs.

\section{Preliminary Results and Definitions}

We first make some observations.

\begin{prop} \label{P:unlinked}
For any $n$, the graphs $K_{n,1}$, $K_{n, 2}$, $K_{n,3}$, $K_{n, 1, 1}$, $K_{n, 2, 1}$ and $K_{n,1,1,1}$ have linkless embeddings.
\end{prop}
\begin{proof}
All of these graphs are subgraphs of $K_{n,1,1,1}$.  However, any cycle in $K_{n,1,1,1}$ must use at least two of the vertices of degree $n+2$, so there are no pairs of disjoint cycles, and hence no links. \end{proof}

The next result is due to the fact that $K_6$ is a minor-minimal intrinsically linked graph \cite{mrs}.

\begin{prop} \label{P:6vertex}
The only intrinsically linked graph with six or fewer vertices is $K_6$, which can be embedded with exactly one non-trivial link.
\end{prop}

\begin{defn} \label{D:numberoflinks}
Given a graph $G$, we define $mnl(G)$ to be the minimal number of links in any embedding of $G$ in $\mathbb{R}^3$.
\end{defn}

Our results are summarized in Table~\ref{Ta:results}.  Since we have only considered graphs with 8 or fewer vertices, all links have two components (three disjoint cycles requires at least 9 vertices).  And since almost all of our arguments are based on the linking number modulo 2, we are really counting the number of two-component links with odd linking number.  For most of these graphs this is sufficient, but it is known that some graphs will always have non-trivial links with even linking number \cite{df}.

We found upper bounds for the minimal number of links by computing this number for specific examples.  This was done using {\it Mathematica}, by modifying a program written by Ramin Naimi \cite{na}.  The program computes the linking number of all pairs of cycles.  To check for non-trivial links with linking number 0, such as the Whitehead link, the program also produced a list of all pairs of cycles with more than 4 crossings, which were checked manually.

\begin{table}
	\begin{center}
	\begin{tabular}{| l | c | c | l | c |}
		\cline{1-2} \cline{4-5}
		\multicolumn{2}{| c |}{\bf Graphs on 7 vertices} & \hspace{.5in} & \multicolumn{2}{| c |}{\bf Graphs on 8 vertices}\\ \cline{1-2} \cline{4-5}
		$G$ & $mnl(G)$ & & $G$ & $mnl(G)$\\ \cline{1-2} \cline{4-5}
		$K_{6,1}$ & 0 & & $K_{7,1}$ & 0 \\ \cline{1-2} \cline{4-5}
		$K_{5,2}$ & 0 & & $K_{6,2}$ & 0\\ \cline{1-2} \cline{4-5}
		$K_{4,3}$ & 0 & & $K_{5,3}$ & 0\\ \cline{1-2} \cline{4-5}
		$K_{5,1,1}$ & 0 & & $K_{4,4}$ & 2\\ \cline{1-2} \cline{4-5}
		$K_{4,2,1}$ & 0 & & $K_{6,1,1}$ & 0\\ \cline{1-2} \cline{4-5}                                                                                                                                                                                                                                                                                                                                                                                                                                            
		$K_{3,3,1}$ & 1 & & $K_{5,2,1}$ & 0\\ \cline{1-2} \cline{4-5}
		$K_{3,2,2}$ & 0 & & $K_{4,3,1}$ & 6\\ \cline{1-2} \cline{4-5}
		$K_{4,1,1,1}$ & 0 & & $K_{4,2,2}$ & 2 \\ \cline{1-2} \cline{4-5}
		$K_{3,2,1,1}$ & 1 & & $K_{3,3,2}$ & 17\\ \cline{1-2} \cline{4-5}
		$K_{2,2,2,1}$ & 0 & & $K_{5,1,1,1}$ & 0 \\ \cline{1-2} \cline{4-5}
		$K_{3,1,1,1,1}$ & 3 & & $K_{4,2,1,1}$ & 6 \\ \cline{1-2} \cline{4-5}
		$K_{2,2,1,1,1}$ & 1 & & $K_{3,3,1,1}$ & 25 \\ \cline{1-2} \cline{4-5}
		$K_{2,1,1,1,1,1}$ & 9 & & $K_{3,2,2,1}$ & 28 \\ \cline{1-2} \cline{4-5}
		$K_7$ & 21 & & $K_{2,2,2,2}$ & 3 \\ \cline{1-2} \cline{4-5}
		\multicolumn{3}{c |}{} & $K_{4,1,1,1,1}$ & 12 \\ \cline{4-5}
		\multicolumn{3}{c |}{} & $K_{3,2,1,1,1}$ & $34 \leq mnl(G) \leq 43$ \\ \cline{4-5}
		\multicolumn{3}{c |}{} & $K_{2,2,2,1,1}$ & $30 \leq mnl(G) \leq 42$ \\ \cline{4-5}
		\multicolumn{3}{c |}{} & $K_{3,1,1,1,1,1}$ & $53 \leq mnl(G) \leq 82 $\\ \cline{4-5}
		\multicolumn{3}{c |}{} & $K_{2,2,1,1,1,1}$ & $54 \leq mnl(G) \leq 94$ \\ \cline{4-5}
		\multicolumn{3}{c |}{} & $K_{2,1,1,1,1,1,1}$ & $111 \leq mnl(G) \leq 172$ \\ \cline{4-5}
		\multicolumn{3}{c |}{} & $K_8$ & $217 \leq mnl(G) \leq 305$ \\ \cline{4-5}
	\end{tabular}
	\end{center}
	\caption{Minimum number of links for complete partite graphs on 7 and 8 vertices} \label{Ta:results}
\end{table}

\section{Counting Links in Complete Graphs with 7 vertices} \label{S:7vertex}

\begin{figure}
    $$\includegraphics{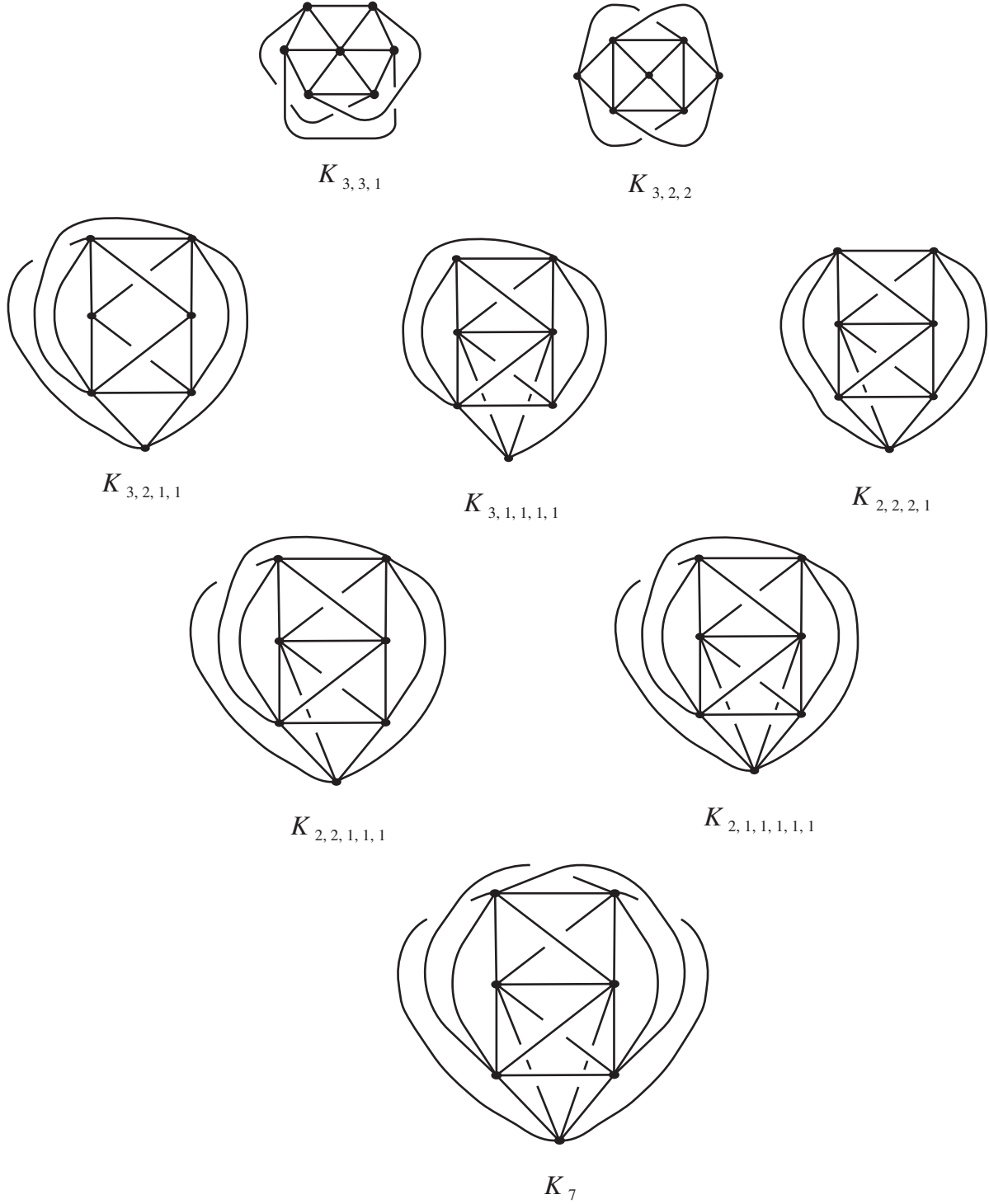}$$
    \caption{Complete partite graphs on 7 vertices} \label{F:7vertex}
    \end{figure}

In this section, we will consider complete $k$-partite graphs on 7 vertices.  By Proposition \ref{P:unlinked}, we need not consider $K_{6,1}$, $K_{5, 2}$, $K_{4,3}$, $K_{5, 1, 1}$, $K_{4, 2, 1}$ or $K_{4,1,1,1}$.  This leaves us with 8 graphs, shown in Figure~\ref{F:7vertex}.

\begin{prop} \label{P:K322,K2221}
The graphs $K_{3,2,2}$ and $K_{2,2,2,1}$ have linkless embeddings.
\end{prop}
\begin{proof}  The embeddings shown in Figure \ref{F:7vertex} are linkless. \end{proof}

\begin{prop} \label{P:K331}
Every embedding of $K_{3,3,1}$ contains a link consisting of a 3-cycle linked with a 4-cycle.  Moreover, $K_{3,3,1}$ can be embedded with exactly one link (so $mnl(K_{3,3,1}) = 1$).
\end{prop}
\begin{proof}  Motwani et. al. \cite{mrs} showed that every embedding of $K_{3,3,1}$ contained a link with odd linking number.  Since every triangle (3-cycle) in $K_{3,3,1}$ must contain the vertex of degree 6, we do not have two disjoint triangles.  So the link must consist of a triangle (containing the preferred vertex) and a square.  The embedding in Figure \ref{F:7vertex} contains exactly one link. \end{proof}

The key idea in the rest our proofs in this section is to look for copies of $K_{3,3,1}$ inside our other graphs.

\begin{prop} \label{P:K3211}
$mnl(K_{3,2,1,1}) = 1$.
\end{prop}
\begin{proof} $K_{3,2,1,1}$ contains $K_{3,3,1}$ as a subgraph, so it must contain at least one link.  The embedding in Figure \ref{F:7vertex} contains exactly one link. \end{proof}

\begin{prop} \label{P:K31111}
$mnl(K_{3,1,1,1,1}) = 3$ 
\end{prop}
\begin{proof} Consider an embedding $F$ of $K_{3,1,1,1,1}$.  Label the vertices of degree six by 1, 2, 3 and 4 (and the other three vertices by 5, 6 and 7).  Then $F$ contains an embedding of $K_{3,3,1}$, using vertex 1 as the preferred vertex.  By Proposition \ref{P:K331}, $F$ contains a link of a triangle (3-cycle) and square (4-cycle) with odd linking number, in which vertex 1 is in the triangle.  At least one of the vertices 2, 3 and 4 is {\it not} in the triangle - without loss of generality, say vertex 4.  Then $F$ contains a second embedding of $K_{3,3,1}$ using vertex 4 as the preferred vertex, so there is a second triangle-square link with odd linking number, this time with vertex 4 in the triangle.  Hence, we have at least two distinct nontrivial triangle-square links in $F$.

Now, let us consider one of these triangle-square links in more detail; the one in which vertex 1 is in the triangle.  If the triangle contains three of the four vertices 1, 2, 3, 4, the remaining four vertices of $K_{3,3,1}$ form a copy of of $K_{3,1}$, which has no cycles.  So the triangle cannot contain three of these vertices, but must contain at least two of them.  Without loss of generality, we can assume that the cycles in the link are 125 and 3647.  Since vertices 3 and 4 are adjacent in $K_{3,1,1,1,1}$, $F$ contains cycles 364 and 473.  We now look at these cycles homologically in the complement of the triangle 125, as elements of $H_1(\mathbb{R}^3 - F(125)) \cong \mathbb{Z}$.  The isomorphism is simply given by the linking number of the cycle with cycle 125.  In homology, $[3647] = [364] + [473]$.  Since $[3647]$ is odd, exactly one of $[364]$ and $[473]$ is also odd, so triangle 125 must link one of the triangles 364 or 473 with odd linking number.  So $F$ contains at least one nontrivial triangle-triangle link.

So $F$ must contain at least 3 non-trivial links.  But the embedding of $K_{3,1,1,1,1}$ in Figure \ref{F:7vertex} has exactly 3 links, so $mnl(K_{3,1,1,1,1}) = 3$. \end{proof}

\begin{prop} \label{P:K22111}
$mnl(K_{2,2,1,1,1}) = 1$
\end{prop}
\begin{proof}  Since $K_{2,2,1,1,1}$ contains $K_{3,3,1}$ as a subgraph, it must contain at least one link.  The embedding shown in Figure \ref{F:7vertex} has exactly one link. \end{proof}

\begin{lem} \label{L:tetrahedron}
Let $G$ be a graph which contains a subgraph $H$ isomorphic to $K_4$ and let $F$ be an embedding of $G$.  If a cycle $C$ in $G$ disjoint from $H$ has odd linking number with a 3-cycle in $H$, then it has odd linking number with 4 cycles in $H$.  Moreover, if $C$ has odd linking number with a 4-cycle $S$ in $H$, then it has odd linking number with two 4-cycles in $H$.
\end{lem}
\begin{proof} 
Consider the subgraph of $F$ induced by $H$.  This subgraph gives a tetrahedron immersed in $\mathbb{R}^3$.  Label the faces of this tetrahedron $T_1, T_2, T_3, T_4$.  Then $[T_1]+[T_2] + [T_3]+[T_4] = 0$ in $H_1(\mathbb{R}^3-[C])$.  An even number of these homology classes must be odd; since we are assuming at least one is odd, either 2 or 4 of them must be odd.  If all 4 are odd, we're done; so say that only $[T_1]$ and $[T_2]$ are odd.  Then the squares $[T_1] + [T_3]$ and $[T_1] + [T_4]$ are distinct 4-cycles with odd linking number with $C$.  So $C$ has odd linking number with 4 cycles in $K_4$.  

Moreover, if $C$ has odd linking number with one 4-cycle, then cannot have odd linking number with all four faces, so by the argument above it will link a second 4-cycle.
\end{proof}

\begin{prop} \label{P:K211111}
$mnl(K_{2,1,1,1,1,1}) = 9$
\end{prop}
\begin{proof} 
Let $F$ be an embedding of $K_{2,1,1,1,1,1}$.  We will show that $F$ must contain at least 3 triangle-triangle links and at least 6 triangle-square links.

Label the two non-adjacent vertices of $K_{2,1,1,1,1,1}$ by $a$ and $b$, and the other vertices 1, 2, 3, 4 and 5.  Then removing either $a$ or $b$ from $F$ leaves us with an embedding of $K_6$, which contains a triangle-triangle link with odd linking number by \cite{cg}.  So we get two distinct nontrivial triangle-triangle links - one link containing $a$ but not $b$, the other containing $b$ but not $a$.  Consider the link containing $a$ - without loss of generality, we may assume the link is between cycles $a12$ and $345$.  Then the subgraph of $F$ induced by the four vertices $b$, 3, 4, and 5 is the embedded 1-skeleton of a tetrahedron.  By Lemma \ref{L:tetrahedron}, $a12$ links 4 cycles in the tetrahedron, including at least one additional triangle.  So there is at least one more nontrivial triangle-triangle link, with $a$ in one triangle and $b$ in the other, which is distinct from the previous two.  So $F$ contains at least three nontrivial triangle-triangle links.

Now we will consider triangle-square links.  Let $M$ be the set of triangles which we know have odd linking number with a square (so initially, $M = \emptyset$). Since $K_{2,1,1,1,1,1}$ contains $K_{3,3,1}$ as a subgraph with vertex $i$ as the preferred vertex (for $i \in S = \{1,2,3,4,5\}$), it contains a triangle-square link with vertex $i$ in the triangle, by Prop~\ref{P:K331}.  Add this triangle to the set $M$.  Continue until all the vertices in $S$ are contained in at least one triangle in $M$ - at this point, $M$ contains at least $\lceil\frac{5}{3}\rceil = 2$ triangles.

Now consider a vertex $i \in S$ which is contained in only one triangle in $M$, say triangle $T = ijk$.  There are two cases to consider.
\begin{enumerate}
	\item ($\{j,k\} \cap \{a,b\} = \emptyset$) In this case, consider a copy of $K_{3,3,1}$ in $F$ in which the vertices are partitioned $(i)(ab*)(jk*)$.  Then there is another triangle-square link, in which $i$ is in a triangle distinct from $T$.  So we can add this triangle to $M$.
	\item ($\{j,k\} \cap \{a,b\} \neq \emptyset$) Without loss of generality, say $j = a$.  Then consider a copy of $K_{3,3,1}$ in $F$ in which the vertices are partitioned $(i)(abk)(***)$.  Again, we have another triangle-square link, in which $i$ is in a triangle distinct from $T$.
\end{enumerate}

So every vertex in $S$ is contained in at least two triangles in $M$, which means that $M$ contains at least $\lceil\frac{10}{3}\rceil = 4$ triangles.  So $F$ contains at least 4 triangle-square links.

The remainder of our proof consists of two cases.
\begin{enumerate}
	\item ($M$ contains a triangle $T$ with vertex $a$)  In this case, we will show that $T$ links {\it two} squares.  Without loss of generality, say that $T = a12$ links square $b345$ with odd linking number.  The vertices $b, 3, 4, 5$ are all adjacent, so the subgraph they induce is isomorphic to $K_4$.  Then, by Lemma \ref{L:tetrahedron}, $T$ links a second square with odd linking number.
	\item ($M$ does not contain a triangle with vertex $a$) We know that $a$ is contained in at least one triangle-triangle link, say in triangle $T = a12$.  As in Lemma~\ref{L:tetrahedron}, this means $[345] + [b43] + [b54] + [b35] = 0$ in $H_1(\mathbb{R}^3 - F(a12))$, and either two or four of the terms are odd.  If two are odd, we can combine an odd and even term to get a square which links $T$ with odd linking number.  If all four are odd, then at least two are the same sign - combining these gives a square which links $T$ with {\it nonzero} even linking number.  In either case, $T$ links a square.
\end{enumerate}

In both of these cases, we get a new triangle-square link with $a$ in the triangle.  Similarly, we can show there will be a new triangle-square link with $b$ in the triangle.  So we have at least 6 triangle-square links.

The embedding of $K_{2,1,1,1,1,1}$ in Figure \ref{F:7vertex} contains exactly 3 triangle-triangle links and 6 triangle-square links, so $mnl(K_{2,1,1,1,1,1}) = 3 + 6 = 9$. 
\end{proof}

Finally, we consider the complete graph $K_7 = K_{1,1,1,1,1,1,1}$.

\begin{thm} \label{T:K7}
$mnl(K_7) = 21$
\end{thm}
\begin{proof}  
Let $F$ be an embedding of $K_7$.  We will show that $F$ contains at least 21 two-component links with odd linking number.  Since $K_7$ contains 7 distinct copies of $K_6$ (by ignoring each vertex in turn), it contains at least 7 different triangle-triangle links (links where both components are 3-cycles) \cite{cg}.

Using an argument similar to Proposition~\ref{P:K211111}, we will show that there are 7 distinct triangles which each have odd linking number with a square (4-cycle).  Lemma \ref{L:tetrahedron} will then imply that there are at least 14 triangle-square links, completing the proof.

Let $M$ be the set of triangles which we have shown to have odd linking number with a square (so, initially, $M = \emptyset$).  If there is a vertex $i$ in $F$ which has not yet been used in a triangle in $M$, then consider a subgraph of $F$ isomorphic to $K_{3,3,1}$ which has $i$ as the preferred vertex (the vertex of degree 6).

By Proposition~\ref{P:K331}, there is a link in our subgraph with odd linking number, consisting a triangle through vertex $i$ and a square.  Add this triangle to $M$.  Since $i$ was not previously used, this triangle will not yet be an element of $M$.

Continue this process until every vertex has been used at least once.  Since $\lceil\frac{7}{3}\rceil = 3$, $M$ will contain at least 3 triangles.  Now consider a vertex $i$ which is used in {\it exactly} one triangle $T = ijk$ in $M$.  Consider a subgraph of $K_{3,3,1}$ in $F$ where the vertices are partitioned $(i)(jk*)(***)$, so the subgraph does not contain the edge $jk$, and so does not contain the triangle $T$.  This subgraph will contain a link with odd linking number, consisting of a triangle through vertex $i$ and a square.  This triangle can be added to $M$, since it is not $T$, which was the only triangle in $M$ containing vertex $i$.

We can continue this process until every vertex is used in at least 2 triangles in $M$.  At this point, $M$ will contain at least $\lceil \frac{14}{3} \rceil = 5$ triangles.  Now suppose that vertex $i$ is used in {\it exactly} two triangles, $T_1$ and $T_2$.  There are two cases, depending on whether $T_1$ and $T_2$ share an edge, or only a vertex.
\begin{enumerate}
	\item $T_1 = ijk$ and $T_2 = ijl$, so the two triangles share an edge.  Then we consider the $K_{3,3,1}$ inside $F$ formed using the partition $(i)(jkl)(***)$, which contains neither $T_1$ nor $T_2$, but will contain a triangle-square link involving a third triangle $T_3$ through vertex $i$.  We can add $T_3$ to $M$.
	\item $T_1 = ijk$ and $T_2 = ilm$, so the two triangles only share the vertex $i$.  Then consider the $K_{3,3,1}$ inside $F$ formed using the partition $(i)(jk*)(lm*)$, which contains neither $T_1$ nor $T_2$.  Then, as in the last case, we will get a new triangle $T_3$ containing $i$ which we can add to $M$.
\end{enumerate}

So, ultimately, every vertex will be used in at least 3 triangles in $M$, so the set will contain at least $\lceil \frac{21}{3} \rceil = 7$ triangles.  So there are at least 7 distinct triangles in $F$ which have odd linking number with a square.  By Lemma \ref{L:tetrahedron}, this gives us at least 14 different triangle-square links.  In addition, there are at least 7 different triangle-triangle links, for a total of 21 distinct links. 

Figure \ref{F:7vertex} shows an embedding of $K_7$ which contains exactly 21 links, which shows that $mnl(K_7) = 21$. 
\end{proof}

\section{Counting Links in Complete Graphs with 8 vertices} \label{S:8vertex}

\begin{figure}
    $$\includegraphics{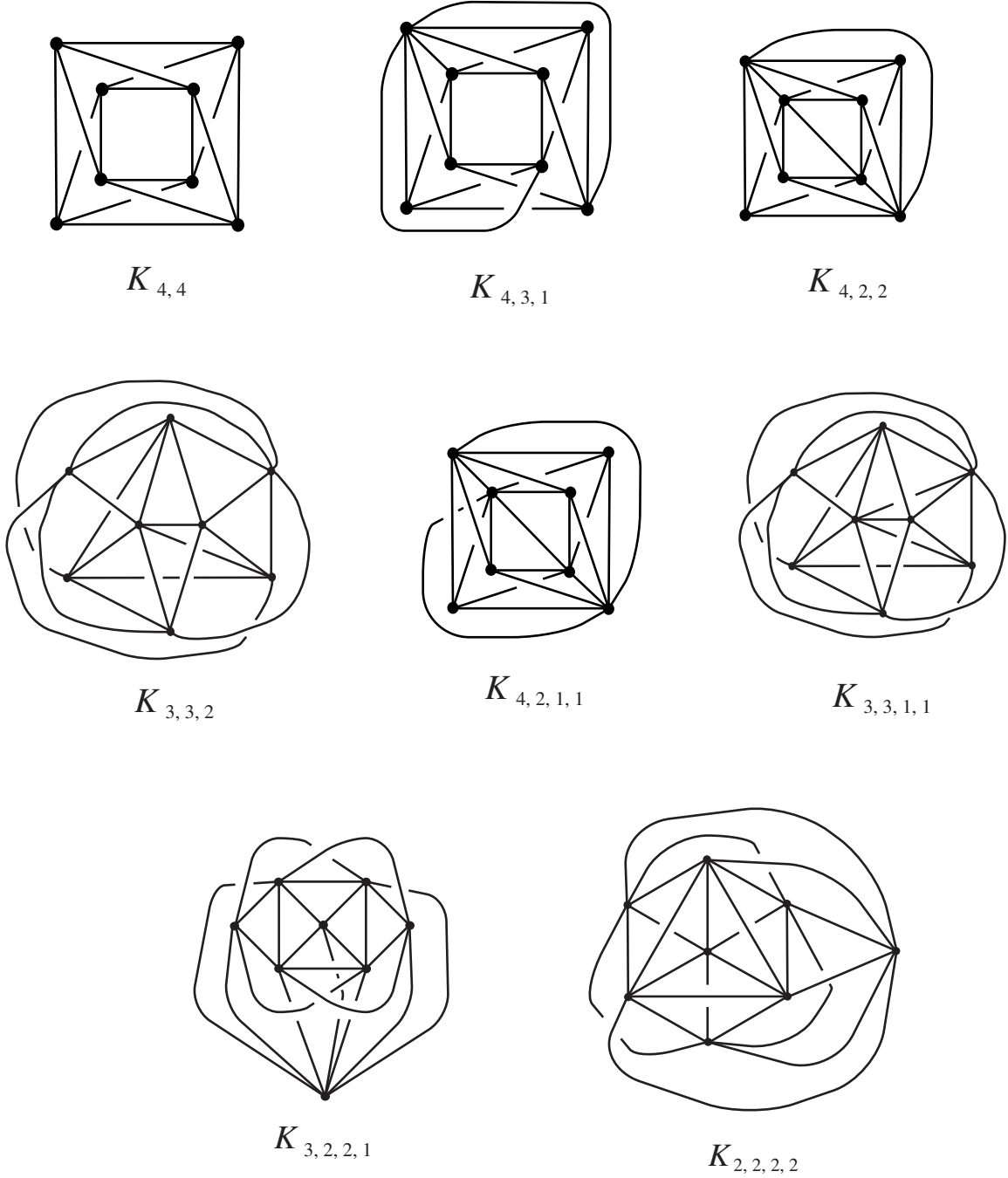}$$
    \caption{Complete partite graphs on 8 vertices} \label{F:8vertex}
    \end{figure}

\begin{figure}
    $$\scalebox{.9}{\includegraphics{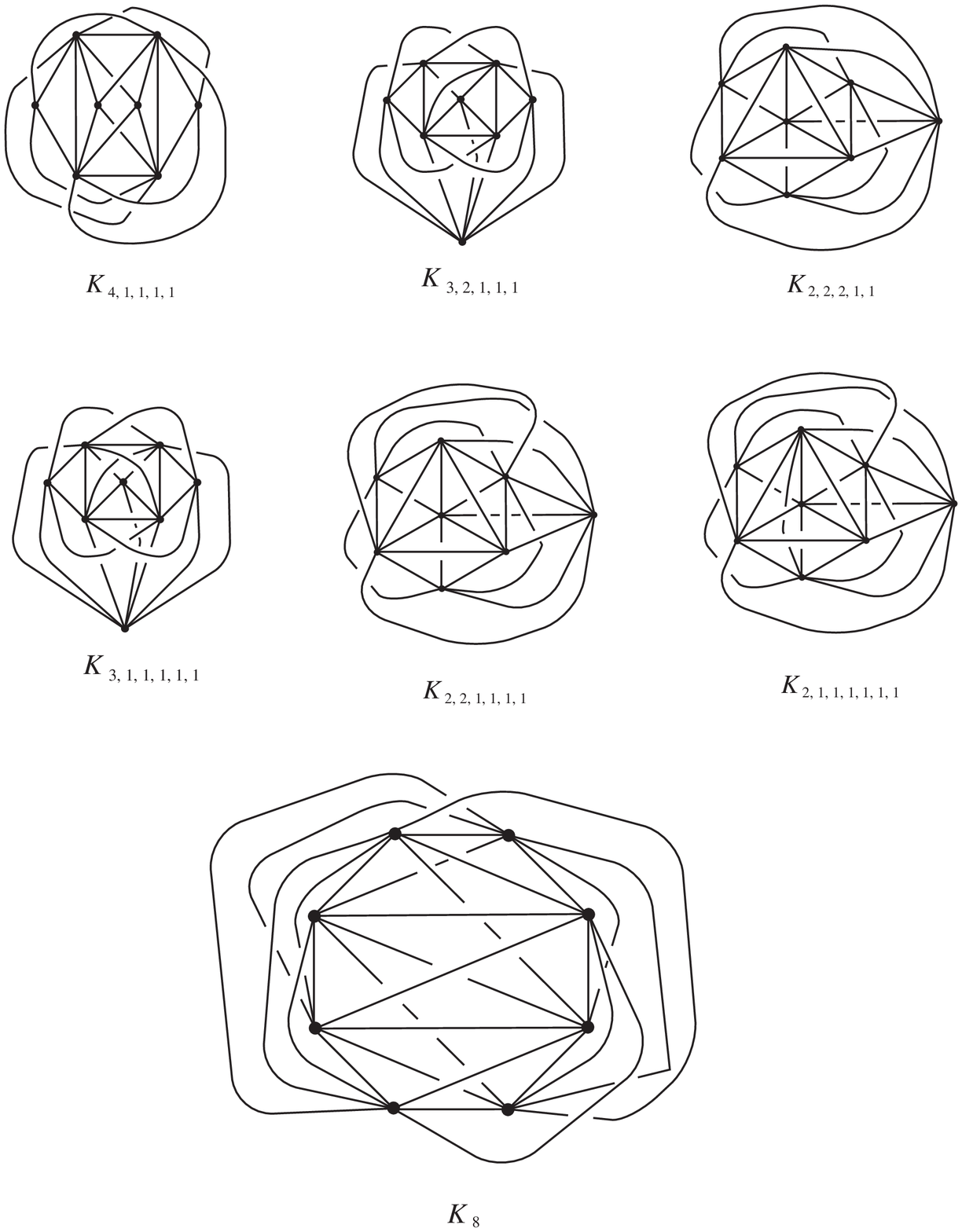}}$$
    \caption{More complete partite graphs on 8 vertices} \label{F:8vertex2}
    \end{figure}

Now we turn to complete graphs on 8 vertices.  By Proposition \ref{P:unlinked}, we do not need to consider $K_{7,1}$, $K_{6, 2}$, $K_{5,3}$, $K_{6, 1, 1}$, $K_{5, 2, 1}$ or $K_{5, 1, 1, 1}$.  This leaves us with 15 other graphs; these are shown in Figures \ref{F:8vertex} and \ref{F:8vertex2}.

The following lemma will be useful in several of our proofs.

\begin{lem} \label{L:K44edge}
Given any embedding $F$ of $K_{4,4}$ and any edge $e$, there is a link in $F$ containing the edge $e$.
\end{lem}
\begin{proof}
Without loss of generality, partition the vertices of $K_{4,4}$ as $(1357)(2468)$, and let $e$ be the edge 78.  If we contract $e$ in $F$ we get an embedding of $K_{3,3,1}$, which contains a link of a triangle (passing through the contracted vertex 7/8) and a square.  This lifts to a link in $F$ of two squares, with $e$ in one of the squares, as desired.
\end{proof}0

\begin{prop} \label{P:K44}
$mnl(K_{4,4}) = 2$
\end{prop}
\begin{proof}
By Lemma \ref{L:K44edge}, every edge of $K_{4,4}$ is in a link.  The only cycles in $K_{4,4}$ are squares, so every link is between 2 squares, and involves 8 edges.  There are a total of 16 edges in $K_{4,4}$, so for every edge to be in a link, we must have at least two different square-square links (in the minimal case, no two of these squares share an edge).  The embedding shown in Figure \ref{F:8vertex} contains exactly two links, so $mnl(K_{4,4}) = 2$. 
\end{proof}

\begin{prop} \label{P:K431}
$mnl(K_{4,3,1}) = 6$
\end{prop}
\begin{proof} $K_{4,3,1}$ contains 4 different subgraphs isomorphic to $K_{3,3,1}$ (by choosing 3 of the 4 vertices in the first partition), so any embedding contains at least 4 different triangle-square links, by Proposition \ref{P:K331}.  Moreover, $K_{4,3,1}$ contains a subgraph isomorphic to $K_{4,4}$, so any embedding contains at least 2 different square-square links, by Proposition \ref{P:K44}.  So any embedding of $K_{4,3,1}$ contains at least 6 links, and the embedding shown in Figure \ref{F:8vertex} has exactly 6. \end{proof}

\begin{prop} \label{P:K422}
$mnl(K_{4,2,2}) = 2$
\end{prop}
\begin{proof}  $K_{4,2,2}$ contains $K_{4,4}$, so $mnl(K_{4,2,2}) \geq 2$, by Proposition \ref{P:K44}.  But the embedding in Figure \ref{F:8vertex} has exactly two links, so $mnl(K_{4,2,2}) = 2$. \end{proof}

\begin{lem} \label{L:pyramid}
Let $G$ be a graph which contains a subgraph $H$ isomorphic to $K_{2,2,1}$ (the 1-skeleton of a pyramid) and let $F$ be an embedding of $G$.  If a cycle $C$ has odd linking number with one of the faces of the pyramid in $F$, then it has odd linking number with at least 6 cycles in the embedding of $H$ in $F$, and links at least two pentagons (possibly one with even linking number).
\end{lem}
\begin{proof}
Say that the vertices of $H$ are $(13)(24)(a)$.  The faces of the pyramid are the cycles $a12$, $a14$, $a23$, $a34$ and $1234$ (see Figure \ref{F:pyramid}).  
    \begin{figure}
    $$\includegraphics{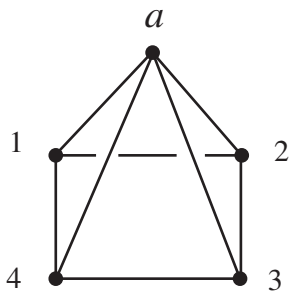}$$
    \caption{The pyramid $K_{2,2,1}$} \label{F:pyramid}
    \end{figure} 
So in $H_1(\mathbb{R}^3 - C)$, the sum $[a12] + [a14] + [a23] + [a34] + [1234] = 0$, which means that an even number of the homology classes are odd.  Since $C$ has odd linking number with at least one face, it must have odd linking number with either 2 or 4 of the faces.  There are several cases to consider.

{\sc Case 1:}  $C$ has odd linking with one triangular face and the square face (without loss of generality, $[a12]$ and $[1234]$).  Then $C$ will also have odd linking with the cycles obtained by adding each of these to each of the other faces (except $[a12]$ and $[a34]$, which are not adjacent), namely $[a123] $, $[a412]$, $[a4123]$, $[a1234]$ and $[a3412]$.  So in this case $C$ links 7 cycles.

{\sc Case 2:} $C$ has odd linking with two adjacent triangular faces (say $[a12]$ and $[a23]$).  Then $C$ also has odd linking with $[a412]$, $[a234]$, $[a2341]$ and $[a3412]$.  So $C$ links 6 cycles.

{\sc Case 3:} $C$ has odd linking with two non-adjacent triangular faces (say $[a12]$ and $[a34]$).  Then $C$ also has odd linking with $[a123]$, $[a412]$, $[a341]$, $[a234]$, $[a2341]$ and $[a4123]$.  So $C$ links 8 cycles.

{\sc Case 4:} $C$ has odd linking with all four triangular faces.  Then $C$ also links the four pentagons formed by adding the base square to each of these triangles.  So $C$ links a total of 8 cycles.

{\sc Case 5:} $C$ has odd linking with three triangular faces and the square face (say $[a12]$, $[a23]$, $[a34]$ and $[1234]$).  Then $C$ also has odd linking with $[a1234]$, $[a412]$ and $[a341]$ (the results of adding $[a14]$ to each of its adjacent faces).  Also, since $[a12]+[a23]+[a34]+[1234] = 0$, and all four linking numbers are odd, the three triangles cannot {\it all} equal $-[1234]$.  So there is another pentagon which links $C$ with non-zero (even) linking number.  So $C$ links 8 cycles in this case.

In every case, $C$ links at least 6 cycles in $K_{2,2,1}$, including at least two pentagons.
\end{proof}

\begin{prop} \label{P:K332}
$mnl(K_{3,3,2}) = 17$
\end{prop}
\begin{proof} First of all, observe that the embedding of $K_{3,3,2}$ in Figure \ref{F:8vertex} has 17 links (1 triangle-triangle, 6 triangle-square, 6 triangle-pentagon and 4 square-square), so we know that $mnl(K_{3,3,2}) \leq 17$.

Let $F$ be an embedding of $K_{3,3,2}$. Say the the vertices of $K_{3,3,2}$ are partitioned $(135)(246)(ab)$.  Then there are two copies of $K_{3,3,1}$ inside $K_{3,3,2}$ - one using vertex $a$, and the other using vertex $b$.  So, by Proposition \ref{P:K331}, there are two triangle-square links in $F$ with odd linking number, each involving one of the vertices $a$ or $b$ (in the triangle) and the six vertices 1, 2, 3, 4, 5, and 6.  Without loss of generality, say that one of these links is between the triangle $T=a12$ and the square $S=3456$.  Johnson and Johnson \cite{jj} showed that $T$ will also link at least one of the four pentagons $b3456$, $3b456$, $34b56$, or $345b6$.  Similarly, the other triangle-square link ($b$ in the triangle) will induce a triangle-pentagon link with $a$ in the pentagon.  Without loss of generality, assume $T$ links the pentagon $P = b3456$.

The subgraph of $K_{3,3,2}$ induced by the vertices of $P$ is the 1-skeleton of a pyramid in $F$, with faces $3456$, $b43$, $b54$, $b65$ and $b36$.  By Lemma \ref{L:pyramid}, $T$ links at least 6 cycles in this pyramid.  Moreover, since $T$ links the square face, we are in either Case 1 or Case 5 of Lemma \ref{L:pyramid}, so $T$ links at least 7 cycles (including $S$ and $P$).  In Case 1, $T$ links one triangle, two additional squares and two additional pentagons.  In Case 5, $T$ links three triangles, two additional squares and one additional pentagon.

So in either case we have at least 5 new links.  Similarly, the triangle-square link coming from $K_{3,3,1}$ with $b$ in the triangle also gives at least 5 new links, and these two sets of links can overlap in at most one triangle-triangle link.  So there are at least $4 + 9 = 13$ links - at least 1 triangle-triangle, 6 triangle-square, 4 triangle-pentagon and 2 others (either triangle-triangle or triangle-pentagon)

In fact, $F$ must contain at least 6 triangle-pentagon links.  If the 2 undetermined links are triangle-triangle links, then there are two additional triangles involved in links.  By Lemma~\ref{L:pyramid}, each of these must link a pentagon (in fact, two) in the complementary pyramid.  So $F$ has at least $4 + 2 = 6$ triangle-pentagon links.

Finally, we consider square-square links.  Since $K_{3,3,2}$ contains a subgraph isomorphic to $K_{4,4}$-edge (by partitioning the vertices $(135a)(246b)$), it contains at least one square-square link, with $a$ and $b$ in different squares.  Then we can get a second subgraph isomorphic to $K_{4,4}$-edge by grouping $a$ with $2, 4, 6$ and $b$ with $1, 3, 5$, giving a different square-square link.  So $F$ contains at least 2 square-square links.

So far, we know that $F$ contains at least 1 triangle-triangle link, at least 6 triangle-square links, at least 6 triangle-pentagon links and at least 2 square-square links, for a total of 15 links.

Once again, let's consider the link of $T = a12$ with $S = 3456$.  Does $3456$ link $a1b2$ as well?  Since $[a1b2] = [a12] + [b21]$ in $H_1(\mathbb{R}^3-S)$, if $[a1b2] = 0$, then $[b21]=-[a12]\neq 0$.  By Lemma \ref{L:pyramid}, $b21$ will also link a triangle.  Then $a12$ and $b21$ are {\it each} involved in a triangle-triangle link, and since they can't link each other, this forces $F$ to have at least two triangle-triangle links.  And, as above, looking at the other triangle in each of these triangle-triangle links will force two new triangle-pentagon links.  This adds at least three links to the 15 we have, for a total of 18, larger than the known minimum.

So, in a minimal case, $3456$ must link $a1b2$.  Similarly, the triangle-square link with $b$ in the triangle (which now cannot be $b12$ with $3456$) will give another new square-square link with $a$ and $b$ in the same square.  So a minimal diagram must have at least 17 links, completing the proof.  
\end{proof}

\begin{prop} \label{P:K4211}
$mnl(K_{4,2,1,1}) = 6$
\end{prop}
\begin{proof}
Since $K_{4,2,1,1}$ contains $K_{4,3,1}$, it contains at least 6 links by Proposition \ref{P:K431}.  The embedding in Figure \ref{F:8vertex} has exactly 6 links. \end{proof}

\begin{prop} \label{P:K3311}
$mnl(K_{3,3,1,1}) = 25$
\end{prop}
\begin{proof}
We first observe that the diagram of $K_{3,3,1,1}$ in Figure \ref{F:8vertex} has 25 links (1 triangle-triangle, 10 triangle-square, 6 triangle-pentagon and 8 square-square), so $mnl(K_{3,3,1,1}) \leq 25$.

We will partition the vertices of $K_{3,3,1,1}$ as $(135)(246)(a)(b)$.  Then, by Proposition \ref{P:K332}, there are at least 17 links, none of which involve the edge $ab$.

As in the proof of Proposition \ref{P:K332}, we assume without loss of generality that we have a link between triangle $a12$ and square $3456$.  We observed in the proof of Proposition \ref{P:K332} that a minimal diagram for $K_{3,3,2}$ must contain a square-square link between cycles $a1b2$ and $3456$.  We will show that if this does {\it not} occur in our diagram of $K_{3,3,1,1}$, then there must be more than 25 links, so the diagram is not minimal.  

Assume that we do {\it not} have a link between squares $a1b2$ and $3456$; at this point, the embedding of $K_{3,3,2}$ has at least 15 links.  Then, as in Proposition \ref{P:K332}, we have a link between cycles $b12$ and $3456$, and each of $a12$ and $b12$ must be involved in a (now distinct) triangle-triangle link with new triangles, so there are at least 16 links.  Say that $a12$ is linked with $bxy$.  If we look at the ``complementary pyramid" to $bxy$, Lemma \ref{L:pyramid} shows that triangle $bxy$ must link at least 5 other cycles in addition to $a12$.  This is also true for the triangle $awz$ linked with $b12$; this adds 10 new links, for a total of 26 in the embedded $K_{3,3,2}$.  But this is larger than the known minimum for $K_{3,3,1,1}$; so this case can be ignored.

So we may now assume we {\it do} have a link between squares $a1b2$ and $3456$.  In $K_{3,3,1,1}$, where we have the edge $ab$, this means that {\it either} $ab1$ or $ab2$ (but not both) have odd linking number with $3456$, bringing our total number of links to 18.

Also, we must {\it not} have a link between $b12$ and $3456$ (since $3456$ already links $a12$), so the triangle-square link in $K_{3,3,2}$ which is induced by $K_{3,3,1}$ as in Proposition \ref{P:K332} involves some other triangle $bxy$ and a new square $S \neq 3456$.  If either $x$ or $y$ are the vertices 1 or 2, then $bxy$ and $a12$ cannot be linked, meaning that they give rise to distinct triangle-triangle links with two other triangles.  As before (when $a1b2$ and $3456$ were not linked), this will lead to more than 25 links in our diagram.  So, without loss of generality, we can assume that our other triangle square link is between $b34$ and $1256$, and that $b34$ and $a12$ are linked.

By  the same argument as before, this gives a square-square link of square $a3b4$ with $1256$, and hence a new triangle-square link with either $ab3$ or $ab4$ linking $1256$.  So our total number of links is now 19.

Without loss of generality, say that $ab1$ has odd linking with $3456$.  Notice that in $H_1(\mathbb{R}^3 - ab1)$, $[3456] = [3452] + [3256]$, so $ab1$ must also have odd linking with either $3452$ or $3256$, giving us a new triangle-square link.  Similarly, we will get another triangle-square link involving $ab3$ or $ab4$.  This gives us a total of 21 links.

So we have at least 4 links where a 3-cycle $abx$ is linked with a 4-cycle $S$; call the remaining vertex $y$.  Then $a$, $b$, $x$ and $y$ form a tetrahedron; by Lemma \ref{L:tetrahedron}, $S$ links 4 cycles in this tetrahedron.  In each case, one of these cycles gives a link we have not previously counted (we leave the details to the reader).  This leaves us with at least 4 new links, for a grand total of 25.  So $mnl(K_{3,3,1,1}) = 25$.
\end{proof}

\begin{lem} \label{L:K2111}
Let $G$ be a graph which contains $K_{2,1,1,1}$ as a subgraph and let $F$ be an embedding of $G$.  If a cycle $C$ has odd linking number with a triangle of $K_{2,1,1,1}$ in $F$, then it has odd linking number with at least 8 cycles in $F$.
\end{lem}
\begin{proof} 
Notice that $K_{2,1,1,1}$ is the 1-skeleton of two tetrahedra joined along one face.  This gives a two-cycle trivial in $H_1(\mathbb{R}^3 - C)$, so the sum of the faces is homologically trivial.  This is shown in Figure \ref{F:2cycle} (here $a$ and $d$ are the vertices of degree 3).  
    \begin{figure}
    $$\includegraphics{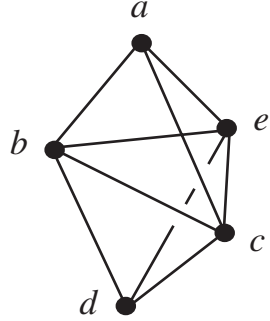}$$
    \caption{$K_{2,1,1,1}$} \label{F:2cycle}
    \end{figure} 
This 2-cycle has six triangular faces: $abc$, $ace$, $aeb$, $dcb$, $dbe$ and $dec$.  Homologically, we have that $[abc] + [ace] + [aeb] + [dcb] + [dbe] + [dec] = 0$.  Since $C$ links at least one of these faces with odd linking number, it must link 2, 4 or 6 of the faces with odd linking number.

There are 9 squares in $K_{2,1,1,1}$, each formed by joining two adjacent faces.  Three of these contain $a$ but not $d$:  $[abce] = [abc] + [ace]$, $[aceb] = [ace] + [aeb]$ and $[aebc] = [aeb] + [abc]$.  Three contain $d$ but not $a$:  $[dbce] = [dbc] + [dce]$, $[dceb] = [dce] + [deb]$ and $[debc] = [deb] + [dbc]$.  Finally, three contain both $a$ and $d$: $[abdc] = [abc] + [bdc]$, $[acde] = [ace] + [cde]$ and $[aedb] = [aeb] + [edb]$.

There are also 6 pentagons in $K_{2,1,1,1}$, each formed by joining three triangular faces:
$$[abcde] = [abc] + [ace] + [dec] = [abe] + [dbc] + [deb]$$
$$[abdce] = [abc] + [dcb] + [ace] = [dce] + [abe] + [deb]$$
$$[abdec] = [abe] + [deb] + [aec] = [abc] + [dec] + [dcb]$$
$$[abedc] = [abe] + [dce] + [aec] = [abc] + [dbe] + [dcb]$$
$$[acbde] = [acb] + [abe] + [deb] = [ace] + [dec] + [dcb]$$
$$[acdbe] = [acb] + [abe] + [dbc] = [ace] + [dec] + [dbe]$$

There are several cases to consider:
\begin{enumerate}
	\item ($C$ links 2 faces sharing a 3-valent vertex)  Without loss of generality, assume that $[abc]$ and $[ace]$ are odd; then $C$ has odd linking with 4 squares and 4 pentagons, for a total of 10 links.
	\item ($C$ links 2 faces {\it not} sharing a 3-valent vertex, but sharing an edge)  Without loss of generality, assume that $[abc]$ and $[bcd]$ are odd; then $C$ has odd linking with 4 squares and 2 pentagons, for a total of 8 links.
	\item ($C$ links 2 faces not sharing an edge)  Without loss of generality, assume that $[abc]$ and $[dbe]$ are odd; then $C$ has odd linking with 6 squares and 4 pentagons, for a total of 12 links.
	\item ($C$ links 4 faces, with three sharing a 3-valent vertex)  Without loss of generality, assume that $[abc]$, $[ace]$, $[abe]$ and $[dec]$ are odd; then $C$ has odd linking with 4 squares and 2 pentagons, for a total of 10 links.
	\item ($C$ links 4 faces, all sharing a 4-valent vertex)  Without loss of generality, assume that $[abc]$, $[ace]$, $[bcd]$ and $[dec]$ are odd (the four faces sharing vertex $c$); then $C$ has odd linking with 4 squares and 4 pentagons, for a total of 12 links.
	\item ($C$ links 4 faces, with two sharing vertex $a$, two sharing vertex $d$, and not all sharing a 4-valent vertex)  Without loss of generality, assume that $[abc]$, $[ace]$, $[dcb]$ and $[dec]$ are odd; then $C$ has odd linking with 6 squares and 2 pentagons, for a total of 12 links.
	\item ($C$ links 6 faces)  In this case, $C$ has odd linking with all the 3-cycles, none of the 4-cycles, and all of the 5-cycles, for a total of 12 links.
\end{enumerate}

So $C$ links at least 8 cycles in $K_{2,2,2,1}$.
\end{proof}

\begin{prop} \label{P:K3221}
$mnl(K_{3,2,2,1}) = 28$
\end{prop}
\begin{proof}
The diagram for $K_{3,2,2,1}$ shown in Figure \ref{F:8vertex} has 28 links (2 triangle-triangle, 10 triangle-square, 10 triangle-pentagon and 6 square-square), so we know that $mnl(K_{3,2,2,1}) \leq 28$.

Assume the vertices of $K_{3,2,2,1}$ are partitioned (123)(45)(67)(8).  First, let's consider the square-square links.  There are several subgraphs isomorphic to $K_{4,4} - {\rm edge}$, which always contains at least one square-square link \cite{mrs}.  First, partition the vertices (1234)(5678) to get a square-square link (where 5 is adjacent to 1, 2 or 3), and then repartition them (1235)(4678) to get another (where 5 is not adjacent to 1, 2 or 3).  In both of these, 6 and 7 are each adjacent to at least one of vertices 1, 2, or 3.  So we can partition the vertices (1236)(4578) and (1237)(4568) to get two new square-square links.  In all of the links we have found so far, vertex 8 was adjacent to at least one of 1, 2, or 3 in one of the squares.  If we consider the partition (1238)(4567) we have a subgraph isomorphic to $K_{4,4}$ in which 8 is not adjacent to 1, 2, or 3.  By Proposition \ref{P:K44} this subgraph contains two new square-square links.  So we have a total of at least 6 square-square links in any embedding of $K_{3,2,2,1}$.

Now we observe that $K_{3,2,2,1}$ contains two subgraphs isomorphic to $K_{3,3,2}$ -- one obtained by partitioning the vertices (123)(458)(67), and the other by partitioning the vertices (123)(678)(45).  Each of these subgraphs contains at least 13 triangle-triangle, triangle-square and triangle-pentagon links (involving two different triangles) by the proof of Proposition \ref{P:K332}; the question is the extent to which these overlap.

We know that the copy of $K_{3,3,2}$ determined by the partition (123)(458)(67) contains two different triangles involved in links (and possibly others); denote these $T_1$ and $T_2$ (vertex 6 is in $T_1$, vertex 7 is in $T_2$).  We first consider the case that neither $T_1$ nor $T_2$ contain vertex 8 (the unique 7-valent vertex in $K_{3,2,2,1}$).  Without loss of generality, say that $T_1 = 146$ and $T_2 = 257$, so $T_1$ links the square $S_1 = 2538$ and $T_2$ links the square $S_2 = 1438$.  But since our graph contains edges 58 and 48, this means that $T_1$ must link either $258$ or $358$, and $T_2$ must link either $148$ or $348$.  Call the new triangles $T_3$ and $T_4$.  Each of these triangles has a complementary pyramid in $K_{3,2,2,1}$ - since each links one face, they must each link at least 5 additional cycles in their respective pyramids (by Lemma \ref{L:pyramid}), giving at least 10 new links.  This means we have at least $6 + 13 + 10 = 29$ links, which is larger than the known minimum of 28.

So either $T_1$ or $T_2$ must contain vertex 8.  Without loss of generality, say $T_1 = 146$ and $T_2 = 287$.  But then $T_2$ is not a cycle in the copy of $K_{3,3,2}$ determined by the partition (123)(678)(45) (although $T_1$ is).  So this graph must contain another triangle $T_3$ passing through vertex 8 (but {\it not} vertex 7) which is involved in links, by the same argument as before.  By the argument in the proof of Proposition \ref{P:K332}, $T_2$ and $T_3$ are each contained in one triangle-triangle link (which must be different, since $T_2$ and $T_3$ are not disjoint cycles), 3 triangle-square links and 3 triangle-pentagon links.  So we have a total of $6 + 7 + 7 = 20$ links.  

It remains to count the links involving triangle $T_1$.   The subgraph induced by the other vertices is isomorphic to $K_{2,1,1,1}$ (i.e. $K_5 - {\rm edge}$).  By Lemma \ref{L:K2111}, $T_1$ links at least 8 cycles in this subgraph.  However, to avoid introducing new triangles (which would force more than 28 links), it will link the faces $T_2$ and $T_3$, which do {\it not} share any of the edges 78, 58 or 57; so if we look at the cases in the proof of Lemma \ref{L:K2111}, we see that $T_1$ must link at least 10 cycles.  At most two of these are the triangle-triangle links we have already counted, leaving us with 8 new links.  This brings our total to 28.  Since we have an example with exactly 28 links, we know that $mnl(K_{3,2,2,1}) = 28$.
\end{proof}

\begin{prop} \label{P:K2222}
$mnl(K_{2,2,2,2}) = 3$
\end{prop}
\begin{proof}  Say the vertices of $K_{2,2,2,2}$ are partitioned (12)(34)(56)(78).  Then there is a subgraph isomorphic to $K_{4,4}$ by grouping the vertices (1234)(5678); by Proposition \ref{P:K44}, this subgraph has at least 2 square-square links.  Without loss of generality, say that one of these links uses edge 15; then we can get a different subgraph isomorphic to $K_{4,4}$ by grouping vertices (1256)(3478), which does not contain the edge 15.  This subgraph {\it also} has at least two square-square links, and at least one of these must be new.  So $K_{2,2,2,2}$ must have at least 3 square-square links.  The embedding in Figure \ref{F:8vertex2} has exactly 3 square-square links, so $mnl(K_{2,2,2,2}) = 3$. \end{proof}

\begin{prop} \label{P:K41111}
$mnl(K_{4,1,1,1,1}) = 12$
\end{prop}
\begin{proof} 
Say that the vertices of $K_{4,1,1,1,1}$ are partitioned (1234)(5)(6)(7)(8).  Then, by grouping the vertices (1234)(5678), we have a subgraph isomorphic to $K_{4,4}$; so $K_{4,1,1,1}$ has at least 2 square-square links by Proposition \ref{P:K44}.

There are also several subgraphs isomorphic to $K_{3,3,1}$.  For example, if we remove vertex 4, we can group the remaining vertices (123)(567)(8); by Proposition \ref{P:K331}, this subgraph will have a triangle-square link.  Without loss of generality, say this link is between cycles 158 and 2637.  But then we can also look at the subgraph induced by the same seven vertices using the partition (123)(568)(7).  This subgraph does not contain the edge 58, so it gives us a second triangle-square link, with a new triangle.  Similarly, we get two new triangle-square links by removing each of vertices 1, 2 and 3, so we have a total of at least 8 triangle-square links.

Finally, consider once again the link between cycles 158 and 2637.  In $K_{4,1,1,1,1}$, we also have the edge 67, so, by our usual homology argument, 158 must link either triangle 267 or 367.  Moreover, if we now look in the subgraph induced by (234)(567)(8), we get another triangle-triangle link which does not involve vertex 1.  So we have at least 2 triangle-triangle links.

So an embedding of $K_{4,1,1,1,1}$ must contain at least $2+8+2=12$ links.  The embedding shown in Figure \ref{F:8vertex2} has exactly 12 links (2 triangle-triangle, 8 triangle-square and 2 square-square), so $mnl(K_{4,1,1,1,1}) = 12$. 
\end{proof}

\begin{prop} \label{P:K32111}
$34 \leq mnl(K_{3,2,1,1,1}) \leq 43$
\end{prop}

\begin{proof}
The embedding shown in Figure \ref{F:8vertex2} has exactly 43 links (3 triangle-triangle links, 16 triangle-square links, 10 square-square links and 14 triangle-pentagon links), so $mnl(K_{3,2,1,1,1}) \leq 43$.

The graph $K_{3,2,1,1,1}$ has $K_{3,2,2,1}$ as a subgraph, and hence
every embedding contains at least 28 non-split links, by Proposition \ref{P:K3221} (2 triangle-triangle, 10 triangle-square, 10 triangle-pentagon and 6 square-square).  Suppose
$K_{3,2,1,1,1}$ is partitioned (123)(45)(6)(7)(8), and the subgraph isomorphic to $K_{3,2,2,1}$
is (123)(45)(67)(8).  Then any link using edge 67 will be new.  We
can find a subgraph of the form $K_{4,4}$ by partitioning
(1236)(4578), and by Lemma \ref{L:K44edge}, there is a square-square link
using edge 67, in which 7 is adjacent to 1, 2 or 3.  Partitioning
(1237)(4568) we find another $K_{4,4}$ containing edge 67, but in
the square-square link that uses this edge, 7 is not adjacent to 1,
2 or 3.  Thus, we have found two new square-square links that appear
in every embedding of $K_{3,2,1,1,1}$ in addition to those arising
from $K_{3,2,2,1}$, so the lower bound for $mnl(K_{3,2,1,1,1})$ is at least 30.  Since $K_{3,2,2,1}$ has at least 6 square-square links, $K_{3,2,1,1,1}$ has at least 8 such links.

Now we consider subgraphs of $K_{3,2,1,1,1}$ isomorphic to $K_{3,3,2}$.  There are four such subgraphs, induced by the partitions (123)(456)(78), (123)(457)(68), (123)(458)(67) and (123)(678)(45).  Each of these subgraphs contains a pair of linked triangles, by Proposition \ref{P:K332}.  However, there is no triangle which appears in all 4 subgraphs, so these links require at least 4 distinct triangles in $K_{3,2,1,1,1}$ (3 triangles would require one of the triangles to appear in all four subgraphs, since there is no set of three mutually disjoint triangles among the four subgraphs).  By Proposition \ref{P:K332}, each of these triangles is involved in at least 7 links, at least three of which are triangle-square links.  So $K_{3,2,1,1,1}$ has at least 12 triangle-square links.  Moreover, each of these triangles either links 3 pentagons, or one pentagon and a new triangle which itself links a pentagon, inducing an additional $3\cdot 4 = 12$ links.  Together with the eight square-square links, this means that an embedding of $K_{3,2,1,1,1}$ has at least $2 + 12 + 12 + 8 = 34$ links.
\end{proof}

\begin{prop} \label{P:K22211}
$30 \leq mnl(K_{2,2,2,1,1}) \leq 42$
\end{prop}
\begin{proof}

Say that the vertices of $K_{2,2,2,1,1}$ are partitioned $(12)(34)(56)(7)(8)$; then there is a subgraph isomorphic to $K_{3,2,2,1}$ obtained from the partition $(127)(34)(56)(8)$.  By Proposition \ref{P:K3221}, we will have at least 28 links in $K_{2,2,2,1,1}$, none of which involve edges $17$ or $27$.

There is also a subgraph isomorphic to $K_{4,4}$, using the partition $(1234)(5678)$.  By Lemma \ref{L:K44edge}, there are square-square links involving edges $17$ and $27$.  If these are different, we have two new links.  If these are the same link, then there is a square-square link between square $172x$ and another square $S$, where $x \in \{3, 4, 5, 6, 8\}$.  Since vertices $7$ and $x$ are adjacent in $K_{2,2,2,1,1}$, $S$ will also link either $17x$ or $27x$, and neither of these triangles are in $K_{3,2,2,1}$.  So, once again, we have at least two links which were not in $K_{3,2,2,1}$.

This gives a lower bound of 30 links in $K_{2,2,2,1,1}$.  The embedding shown in Figure \ref{F:8vertex2} has 42 links (3 triangle-triangle, 15 triangle-square, 15 triangle-pentagon and 9 square-square).
\end{proof}

\begin{prop} \label{P:K311111}
$53 \leq mnl(K_{3,1,1,1,1,1}) \leq 82$
\end{prop}

\begin{proof}

Partition the vertices of $K_{3,1,1,1,1,1}$ as (123)(4)(5)(6)(7)(8).
There are three subgraphs of the form $K_{2,1,1,1,1,1}$, obtained by
deleting one of $\{1,2,3\}$.  Each of these subgraphs contains at
least 3 triangle-triangle links. If such a link is contained in two
of these subgraphs, then it cannot be contained in the third.  Thus
we have at least $\lceil \frac{9}{2} \rceil = 5$ triangle-triangle
links.

To count triangle-square links, we delete vertices one at a time.
Deleting one of vertices 1 through 3 leaves a copy of
$K_{2,1,1,1,1,1}$, which contains at least 6 triangle-square links.
Since there are three such graphs, this give 18 triangle-square
links. Deleting one of vertices 4 through 8 leaves a copy of
$K_{3,1,1,1,1}$ which contains at least 2 distinct triangle-square
links, for a total of 10.  Thus, every embedding of
$K_{3,1,1,1,1,1}$ contains at least 10+18 = 28 triangle-square
links.

There are 5 distinct ways to form a $K_{4,4}$ subgraph of
$K_{3,1,1,1,1,1}$, these are of the form $(123i)(****)$.  Each copy of
$K_{4,4}$ contains at least two square-square links.  Notice that as
vertices 1 through 3 must be in the same partition, each of these
square-square links is contained in only one of the $K_{4,4}$
subgraphs above, so we have 10 distinct square-square links.

We now examine subgraphs of the form $K_{3,3,2}$.  We may choose to
partition the vertices such that any vertex from $\{4 \dots 8\}$ is
in the partition of size 2.  Thus, these vertices are contained in
at least one triangle that is used in a triangle-pentagon link, by Proposition \ref{P:K332}.
Suppose it is $ijk$.  Then partition the vertices $(ij)(***)(123)$.
This contains a new triangle-pentagon link, with $i$ used in
the triangle.  So vertices 4 through 8 are used in at least two such
triangles.  However, one vertex of each triangle is taken from the
partition (123).  Thus there are $\frac{10}{2} = 5$ such triangles
involved in triangle-pentagon links.  Since one of $\{1,2,3\}$ is
used in the triangle, the vertices of the pentagon form a $K_{5}
\setminus e$.  Thus each of the triangles links at least 2
pentagons (by Lemma \ref{L:K2111}), for a total of 10 distinct triangle-pentagon links.

This gives a lower bound of 5+28+10+10 = 53 links.  The embedding shown in Figure \ref{F:8vertex2} has 82 links (10 triangle-triangle, 34 triangle-square, 24 triangle-pentagon and 14 square-square), so $53 \leq mnl(K_{3,1,1,1,1,1}) \leq 82$.
\end{proof}

\begin{lem} \label{L:3pentagon}
Let $G$ be a graph which contains $K_5$ as a subgraph, let $F$ be an embedding of $G$, and let $P$ be a 5-cycle in this $K_5$.  If a cycle $C$ has odd linking number with $P$ in $F$, then it has odd linking number with at least three 5-cycles in $F$.
\end{lem}
\begin{proof} 
Notice that $K_5 - {\rm edge} = K_{2,1,1,1}$.  Say that $P = abcde$.  Consider the $K_{2,1,1,1}$ formed by removing edge $ad$ from the $K_5$ induced by the vertices of $P$, as in Figure \ref{F:2cycle}.  This gives the homology element $[P] = [abcde] = [abc] + [ace] + [dec] = [abe] + [dbc] + [deb]$.  Since $[P]$ is odd, so are an odd number of {\it each} set of three triangles, so $C$ links at least one of the faces.

By Lemma \ref{L:K2111}, this means $C$ links at least 8 cycles in $K_{2,1,1,1}$, including at least one other pentagon.  This pentagon shares at most 3 edges with $P$, and so has at least two edges which $P$ does not.  We can form a new $K_{2,1,1,1}$ by restoring edge $ad$ and removing one of these two edges.  $C$ will then link another pentagon on the boundary of this new 2-cycle, distinct from the other two.  So $C$ links at least 3 pentagons. 
\end{proof}

\begin{prop} \label{P:K221111}
$54 \leq mnl(K_{2,2,1,1,1,1}) \leq 94$
\end{prop}

\begin{proof}
Say the vertices are partitioned (12)(34)(5)(6)(7)(8).  An embedding
of $K_{2,2,1,1,1,1}$ contains 4 distinct copies of $K_{6}$. Thus,
the embedding contains at least 4 triangle-triangle links.

Eliminating one of the vertices 1 through 4 leaves a subgraph of the
form $K_{2,1,1,1,1,1}$.  Each of these contains 6 distinct
triangle-square links by Proposition \ref{P:K211111}.  Eliminating one of vertices 5 through 8 leaves
a graph of the form $K_{2,2,1,1,1}$, each of which contains at least one
triangle-square link by Proposition \ref{P:K22111}.  This gives a minimum of 28 distinct
square-triangle links in any embedding of $K_{2,2,1,1,1,1}$.

By deleting the edges between vertex 8 and vertices 1 and 2, we have
$K_{3,2,1,1,1}$ as a subgraph of $K_{2,2,1,1,1,1}$. There are at
least 8 square-square links in every embedding of $K_{3,2,1,1,1}$ by Proposition \ref{P:K32111}.
In addition, any square-square links that use edges 18 or 28 will be
distinct from these. Form a $K_{4,4}$ by the partition (1234)(5678).
There is a square-square link in this graph that uses edge 18.  If
it does not also use edge 28, we have two new square-square links.
So, suppose that this link does use edge 28. Then without loss of
generality, the link is 1825-3647. Form the partition (1256)(3478).
This copy of $K_{4,4}$ contains edge 18 but does not contain the
square 1825, so there must be some other square-square link using
this edge.  This also gives two new square-square links for a total
of at least 10 in every embedding.

We count triangle-pentagon links by studying subgraphs of the
form $K_{3,3,2}$.  Any vertex can be placed in the partition of size
2, so every vertex is contained in a triangle that is part of such a
link.  Suppose vertex 1 is contained in triangle $1jk$. Then form
the subgraph $(12)(jk*)(***)$ to get a second triangle-pentagon link.  A similar argument can be made for each $i \in \{1,2,3,4\}$, thus each of these vertices is contained in at least
2 distinct triangles, while vertices 5, 6, 7 and 8 are each contained in at least one. This gives a total of $8+4=12$ triangles (not all distinct).  So, we have at least $\lceil \frac{12}{3} \rceil =4$ distinct triangles used in triangle-pentagon links. If there are only four triangles, then each triangle $T$ contains one
vertex from (12) and one from (34).  Thus, the complement of these triangles is $K_{5}$, so each $T$ links at least 3 pentagons by Lemma \ref{L:3pentagon}.  If some triangle contains only one of $\{1,2,3,4\}$, then the complement of that triangle is $K_{5} \setminus e$, so it links with only 2
pentagons, but then, as each of $\{1,2,3,4\}$ is contained in two triangles, there is at least one more triangle-pentagon link. In either case, we have at least 12 triangle-pentagon links.

This gives a lower bound of 4+28+10+12=54 links.  The embedding shown in Figure \ref{F:8vertex2} has 94 links (8 triangle-triangle, 34 triangle-square, 34 triangle-pentagon and 18 square-square), so $54 \leq mnl(K_{3,1,1,1,1,1}) \leq 94$.
\end{proof}

\begin{prop} \label{P:K2111111}
$111 \leq mnl(K_{2,1,1,1,1,1,1}) \leq 172$
\end{prop}

\begin{proof}

Say the vertices of $K_{2,1,1,1,1,1,1}$ are partitioned
(12)(3)(4)(5)(6)(7)(8).  Given an embedding of $K_{2,1,1,1,1,1,1}$,
vertices 1 through 7 form a copy of $K_{2,1,1,1,1,1}$, which
contains 3 triangle-triangle links by Proposition \ref{P:K211111}.  Clearly,
these triangles use vertices 1 through 7.  There are 10 $K_{6}$
subgraphs that contain vertex 8, formed by omitting one of $\{1,2\}$
and one of $\{3 \ldots 7\}$, and one formed by vertices 3 through 8.
Each of these contains a triangle-triangle link using vertex 8, so
there are at least 3+11=14 triangle-triangle links.

Omitting vertices one at time, we obtain two distinct copies of
$K_{7}$, and six copies of $K_{2,1,1,1,1,1}$.  By Theorem
\ref{T:K7} and Proposition \ref{P:K211111}, each $K_{7}$ contains 14 distinct triangle-square
links, and $K_{2,1,1,1,1,1}$ contains 6 such links.  Thus the
embedding of $K_{2,1,1,1,1,1,1}$ has at least 64 triangle-square
links.

To find subgraphs of the form $K_{4,4}$, we will choose two of $\{3
\ldots 8 \}$ and group them with (12).  This gives ${6 \choose 2} =
15$ copies of $K_{4,4}$.  Each copy of $K_{4,4}$ contains two
square-square links, but once again, each such link could be
contained in two distinct $K_{4,4}$'s.  Thus, we have at least 15
distinct square-square links.

Look at graphs of the form $K_{3,3,2}$.  We may choose partitions so
that every vertex is in the partition of size two for some subgraph.
Thus every vertex is in a triangle that is used in a
triangle-pentagon link.  Suppose vertex $i$ is in such a triangle,
and that triangle is $ijk$.  Then we can look at the partition $(i*)(jk*)(***)$ to get a second triangle-pentagon link with $i$ in the triangle.  If $i = 1, 2$, then we can get a third triangle as follows.  If the first two triangles containing vertex 1 are of the form $1jk$ and $1kl$, form
$(12)(jkl)(***)$.  If they are of the form $1jk$ and $1lm$ form
$(12)(jk*)(lm*)$.  Thus, vertices 1 and 2 must be contained in at
least three distinct triangles that are used in triangle-pentagon
links, and these triangles contain only 1 or only 2.  Now
we have $\frac{18}{3}  = 6$ distinct triangles used in
triangle-pentagon links, and furthermore each of them must contain
either vertex 1 or vertex 2.  Thus, the complement of one of these triangles
in $K_{2,1,1,1,1,1,1}$ is $K_{5}$, so by Lemma \ref{L:3pentagon}, each triangle
links at least 3 pentagons.  This gives a total of 18
triangle-pentagon links.

This gives us a lower bound of 14 + 64 + 15 + 18 = 111 links.  The embedding shown in Figure \ref{F:8vertex2} has 172 links (16 triangle-triangle, 64 triangle-square, 60 triangle-pentagon and 32 square-square), so $111 \leq mnl(K_{3,1,1,1,1,1}) \leq 172$.
\end{proof}

\begin{thm} \label{T:K8}
$217 \leq mnl(K_8) \leq 305$
\end{thm}
\begin{proof} $K_8$ contains $\binom{8}{6} = 28$ different copies of $K_6$, so at least 28 different triangle-triangle links.  It also contains $\binom{8}{7} = 8$ different copies of $K_7$, so at least $8\cdot 14 = 112$ different triangle-square links (by Theorem \ref{T:K7}).

To find square-square links, we look at copies of $K_{4,4}$ in $K_8$.  There are $\frac{1}{2}\binom{8}{4} = 35$ ways to partition the 8 vertices of $K_8$ into two sets of 4, so there are 35 distinct copies of $K_{4,4}$ in $K_8$.  Each copy of $K_{4,4}$ will contribute at least two square-square links to an embedding of $K_8$, by Proposition \ref{P:K44}.  However, each link will be contained in two different copies of $K_{4,4}$ - for example, a link between squares 1234 and 5678 would be in the $K_{4,4}$'s arising from the partition $\{1,3,5,7\}\{2,4,6,8\}$ and from the partition $\{1,3,6,8\}\{2,4,5,7\}$.  So there will be at least 35 different square-square links in $K_8$.

To count triangle-pentagon links, we look at copies of $K_{3,3,2}$ inside $K_8$.  By Proposition \ref{P:K332}, an embedding of $K_{3,3,2}$, where $a$ and $b$ are the vertices of degree 6, contains at least one triangle-pentagon link with $a$ in the triangle and $b$ in the pentagon, and another with $b$ in the triangle and $a$ in the pentagon.  Let $M$ be the set of triangles we know are involved in triangle-pentagon links (so initially $M = \emptyset$).  Let $i$ be a vertex which does not appear in a triangle of $M$.  Then we can choose a subgraph of $K_8$ isomorphic to $K_{3,3,2}$ in which $i$ has degree 6, so there is a triangle containing $i$ which is part of a triangle-pentagon link.  Continuing until every vertex is used, we get at least $\lceil\frac{8}{3}\rceil = 3$ triangles.

These triangles have a total of 9 vertices, so some vertices are used only once.  Say vertex $i$ is only in triangle $ijk$.  Then we can choose a partition $(i*)(jk*)(***)$ of the vertices of $K_8$ to get another copy of $K_{3,3,2}$ which does not contain the edge $jk$.  So $i$ will be in a different triangle in another triangle-pentagon link, and we can add this new triangle to $M$.  We can continue in this way until every vertex is used at least twice, giving at least $\lceil\frac{16}{3}\rceil = 6$ triangles.

These triangles have a total of 18 vertices, so some vertices are used only twice.  Say that $i$ is in triangles $T_1 = ijk$ and $T_2 = ilm$.  There are two cases, depending on whether $T_1$ and $T_2$ share an edge, or only a vertex.
\begin{enumerate}
	\item $T_1 = ijk$ and $T_2 = ijl$, so the two triangles share an edge.  Then we consider the $K_{3,3,2}$ inside $K_8$ formed using the partition $(i*)(jkl)(***)$, which contains neither $T_1$ nor $T_2$, but will contain a triangle-pentagon link involving a third triangle $T_3$ through vertex $i$.  We can add $T_3$ to $M$.
	\item $T_1 = ijk$ and $T_2 = ilm$, so the two triangles only share the vertex $i$.  Then consider the $K_{3,3,2}$ inside $K_8$ formed using the partition $(i*)(jk*)(lm*)$, which contains neither $T_1$ nor $T_2$.  Then, as in the last case, we will get a new triangle $T_3$ containing $i$ which we can add to $M$.
\end{enumerate}
Continuing until every vertex is used three times, we get at least $\lceil\frac{24}{3}\rceil = 8$ triangles in $M$.

These triangles have a total of 24 vertices, so some vertices are used only three times.  Say vertex $i$ is in triangles $T_1 = ijk$, $T_2 = ilm$ and $T_3 = ino$.  We have several cases.
\begin{enumerate}
	\item $T_1 = ijk$, $T_2 = ijm$ and $T_3 = ijo$, so all three triangles share an edge $ij$.  Then use the partition $(ij)(***)(***)$ to form a new $K_{3,3,2}$.  In this subgraph, there will be a triangle-pentagon link where $i$ is in the triangle and $j$ is in the pentagon, so the link does not contain edge $ij$.  This gives us a new triangle for $M$.
	\item $T_1 = ijk$, $T_2 = ijm$ and $T_3 = imk$, so each pair of triangles shares an edge.  Then use the partition $(i*)(jmk)(***)$.
	\item $T_1 = ijk$, $T_2 = ijm$ and $T_3 = imo$, so two pairs of triangles share an edge.  Then use the partition $(ij)(mok)(***)$.
	\item $T_1 = ijk$, $T_2 = ijm$ and $T_3 = ino$, so just one pair of triangles share an edge.  Then use the partition $(ij)(mk*)(no*)$.
	\item $T_1 = ijk$, $T_2 = ilm$ and $T_3 = ino$, so no triangles share an edge.  Then use the partition $(ij)(lm*)(no*)$.
\end{enumerate}
So we can continue until every vertex is used 4 times, yielding at least $\lceil\frac{32}{3}\rceil = 11$ triangles in $M$.

These triangles have a total of 33 vertices, so some vertices are used only 4 times.  Say vertex $i$ is in triangles $T_1 = ijk$, $T_2 = ilm$, $T_3 = ino$ and $T_4 = ipq$.  The vertices $i,j,k,l,m,n,o,p,q$ cannot all be distinct; the ``worst case" is when only two are the same, say $j$ and $l$.  In this case, use the partition $(ij)(no*)(pq*)$ to get a new triangle with vertex $i$.  Any other case can be dealt with more easily.  So we can continue until every vertex is used 5 times, yielding at least $\lceil\frac{40}{3}\rceil = 14$ triangles in $M$.

Now we need to determine how many pentagons each of these triangles must link.  Consider a triangle-pentagon link with triangle $T$ and pentagon $P$.  The vertices of $P$ induce a subgraph of $K_8$ isomorphic to $K_5$; so by Lemma \ref{L:3pentagon} $T$ links at least 3 pentagons (including $P$).  So there are at least $14 \cdot 3 = 42$ triangle-pentagon links.

Adding this up, we have at least $28 + 112 + 35 + 42 = 217$ links in $K_8$.  The example shown in Figure \ref{F:8vertex2} (motivated by a minimal crossing diagram found by Guy \cite{gu}) has 305 links (28 triangle-triangle, 112 triangle-square, 56 square-square and 109 triangle-pentagon).  So $217 \leq mnl(K_8) \leq 305$. \end{proof}

\section{Counting Links in Complete Bipartite Graphs} \label{S:bipartite}

In this section we will consider complete {\it bipartite} graphs.  In this case, there is a relatively natural spatial embedding of the graph which we conjecture gives the minimal number of links.  We prove this result for the graphs $K_{4,n}$.

\begin{figure}
    $$\scalebox{.8}{\includegraphics{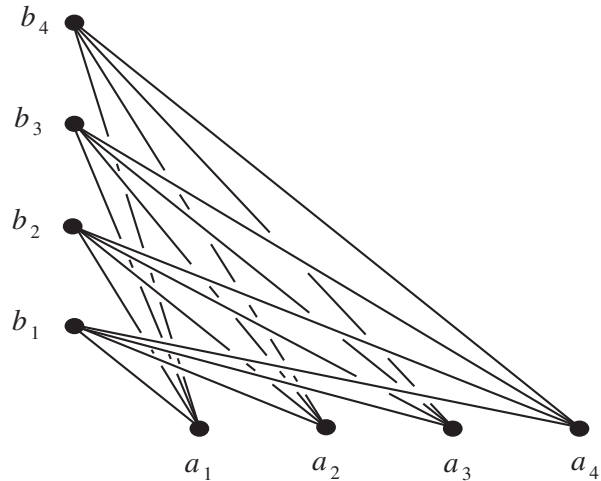}}$$
    \caption{Fan embedding of $K_{4,4}$} \label{F:bipartitefan}
    \end{figure}

We will call this the {\it fan} embedding.  As an example, the fan embedding of $K_{4,4}$ is shown in Figure \ref{F:bipartitefan}.

For a general complete bipartite graph $K_{m,n}$, we can describe the embedding by describing its projection in the plane, as follows.  Denote the two sets of independent vertices by $a_1, \dots, a_m$ and $b_1, \dots, b_n$.  Place the vertex $a_i$ along the $x$-axis at $(i, 0)$ and the vertex $b_j$ along the $y$-axis at $(0, j)$.  Draw the line segments $\overline{a_ib_j}$.  Consider two segments $\overline{a_ib_j}$ and $\overline{a_kb_l}$, where $k > i$.  Then these segments cross if and only if $j < l$, and in this case $\overline{a_ib_j}$ crosses {\it over} $\overline{a_kb_l}$.

\begin{thm} \label{T:K4n}
$mnl(K_{4,n}) = 2 \binom{n}{4}$, and the minimum is realized by the fan embedding.
\end{thm}
\begin{proof}
Label the two sets of independent vertices in $K_{4,n}$ by $a_1, a_2, a_3, a_4$ and $b_1,\dots,b_n$.  The four vertices $a_1, a_2, a_3, a_4$ together with any subset $b_{j_1}, b_{j_2}, b_{j_3}, b_{j_4}$ of the $b_j$'s induce a subgraph of $K_{4,n}$ isomorphic to $K_{4,4}$.  By Proposition \ref{P:K44}, this subgraph contains at least 2 square-square links.  Since each of these links use all 8 vertices in the subgraph, a different subgraph will give different links.  There are $\binom{n}{4}$ ways to choose the vertices $b_{j_1}, b_{j_2}, b_{j_3}, b_{j_4}$, so $K_{4,n}$ must contain at least $2 \binom{n}{4}$ links.

However, in the fan embedding of $K_{4,n}$, the embedding of any such subgraph is isotopic to the fan embedding of $K_{4,4}$ in Figure \ref{F:bipartitefan}, which is isotopic to the embedding of $K_{4,4}$ in Figure \ref{F:8vertex}, and hence has exactly 2 links.  Moreover, any link in $K_{4,n}$ must involve all the $a_i$'s and four of the $b_j$'s, and so is contained in one of these subgraphs.  Therefore, the fan embedding has exactly $2 \binom{n}{4}$ links.  So $mnl(K_{4,n}) = 2\binom{n}{4}$, realized by the fan embedding.
\end{proof}

This result inspires the following conjecture:

\begin{conj}
The fan embedding of $K_{m,n}$ realizes $mnl(K_{m,n})$.
\end{conj}

However, it is not so easy even to compute the number of links in the fan embedding.  It clearly suffices to consider the case of $K_{n,n}$, since any link in $K_{m,n}$, $m > n$, is contained in a subgraph isomorphic to $K_{n,n}$.

\begin{quest}
How many links are in the fan embedding of $K_{n,n}$?
\end{quest}

We know that there are two links in the fan embedding of $K_{4,4}$, and a computer calculation shows that there are 150 links in the fan embedding of $K_{5,5}$ (50 square-square links and 100 square-hexagon links).  The number of links is increasing very rapidly, so computer calculations quickly become infeasible.

\section{Minimal Book Embeddings}

We define an {\it n-book} $B_n$ as the subset of $\mathbb{R}^3$ consisting of a line $L$ (the spine) and $n$ distinct half-planes $S_1$, $S_2$, \dots, $S_n$ (the pages) with $L$ as their common boundary.  

\begin{defn}
Let $G$ be a graph.  An {\bf n-book embedding} of $G$ is a tame embedding $f: G \rightarrow B_n \subset \mathbb{R}^3$ such that:
\begin{enumerate}
	\item $f(V(G)) \subset L$
	\item For each edge $e \in E(G)$, there is exactly one sheet $S_i$ such that $f(e) \subset S_i$.
\end{enumerate}
\end{defn}

A {\it minimal book embedding} of a graph $G$ is a book embedding which minimizes the the number of pages; the {\it pagenumber} of $G$ is the number of pages in a minimal book embedding.  Book embeddings, and particularly minimal book embeddings, minimize the entanglement among the edges of the graph (for example, a book embedding of a graph cannot contain any local knots along the edges).  So it is reasonable to think that minimal book embeddings will also minimize the linking or knotting in an embedding.

Otsuki \cite{ot} gave the first results along these lines.  Since Conway and Gordon \cite{cg} showed that every embedding of $K_6$ contains a pair of linked triangles, and every embedding of $K_7$ contains a knotted 7-cycle, it is immediate that every embedding of $K_n$ contains at least $\binom{n}{6}$ pairs of linked triangles, and $\binom{n}{7}$ knotted 7-cycles.  Otsuki constructed a particular minimal book embedding of $K_n$ called the {\it canonical book representation}, which has the property that removing any vertex gives a canonical book representation of $K_{n-1}$.  He used this property to show that the canonical book representation of $K_n$ contained {\it exactly} $\binom{n}{6}$ pairs of linked triangles and $\binom{n}{7}$ knotted 7-cycles, attaining the minimum possible.

In fact, the embedding of $K_7$ shown in Figure \ref{F:7vertex} is a canonical book representation of $K_7$ (\cite{ot}, Lemma 3.1).  So we can use Theorem \ref{T:K7} to extend Otsuki's result:

\begin{cor} \label{C:K7trisquare}
A canonical book representation of $K_n$ contains exactly $14\binom{n}{7}$ triangle-square links, attaining the minimum possible.
\end{cor}

\begin{proof}
Every $K_7$-subgraph of a canonical book representation of $K_n$ is a canonical book representation of $K_7$, and is therefore ambient isotopic to the embedding of $K_7$ shown in Figure \ref{F:7vertex} (by \cite{ot}, Theorem 1.2).  Therefore, each such subgraph contains 14 triangle-square links, and the embedding of $K_n$ contains $14\binom{n}{7}$ triangle-square links.  By Theorem \ref{T:K7}, this is minimal.
\end{proof}

    \begin{figure}
    $$\scalebox{.8}{\includegraphics{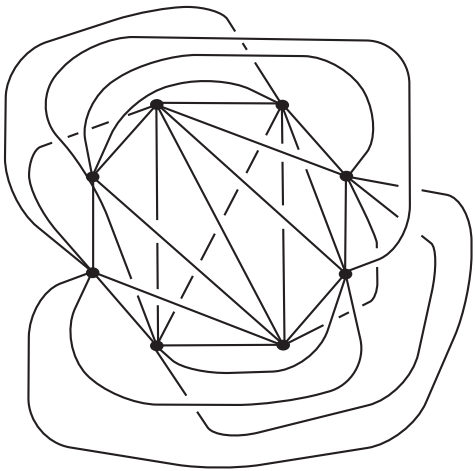}}$$
    \caption{Canonical book representation for $K_8$} \label{F:K8canonical}
    \end{figure} 

We can also look at a canonical book representation for $K_8$, as shown in Figure~\ref{F:K8canonical}.
We can compute that this embedding contains exactly 305 links -- 28 triangle-triangle links, 112 triangle-square links, 112 triangle-pentagon links and 53 square-square links.  This is the same total number of links as the embedding shown in Figure \ref{F:8vertex2}, which leads us to conjecture:

\begin{conj} \label{C:minimal}
For any graph $G$, there is a minimal book embedding of $G$ which realizes $mnl(G)$.
\end{conj}

However, notice that the embeddings of $K_8$ in Figures \ref{F:8vertex2} and \ref{F:K8canonical}, while they have the same total number of links, do {\it not} have the same number of links of each type; so an embedding which minimizes the total number of links may not minimize the number of links of each type.  

\begin{quest}
Given integers $k$ and $l$, is there an embedding of $K_n$ ($n \geq k+l$) which minimizes both the total number of links and the number of links between a $k$-cycle and an $l$-cycle?
\end{quest}

Our discussion of complete bipartite graphs in Section \ref{S:bipartite} also provides some evidence in favor of Conjecture \ref{C:minimal}.  The fan embeddings discussed there are easily seen to be book embeddings, although not minimal (the fan embedding of $K_{n,n}$ has a pagenumber of $n$, while Enomoto et al. \cite{eno} have shown that the pagenumber is at most $\lceil 2n/3 \rceil + 1$).  So once again, it seems reasonable to focus on book embeddings, and if possible on minimal book embeddings, when trying to minimize the total number of links in a graph embedding.

\begin{figure}
    $$\scalebox{.6}{\includegraphics{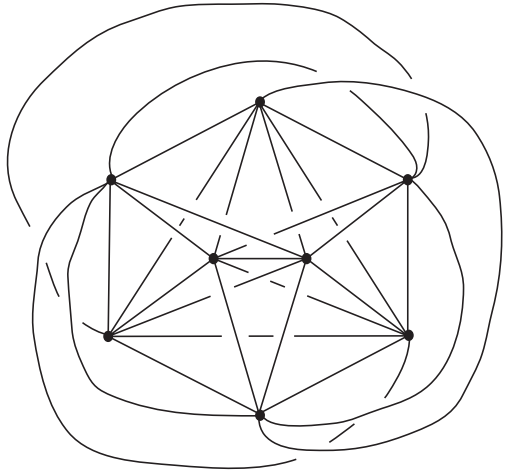}}$$
    \caption{An embedding of $K_8$ with only Hopf links.} \label{F:K8hopf}
    \end{figure}

There are also issues about the interplay between different measures of the complexity of the linking in a graph embedding.  For example, while the embeddings of $K_8$ in Figures \ref{F:8vertex2} and \ref{F:K8canonical} each have 305 links (the smallest number we have found), they also each contain a link with linking number 2.  On the other hand, there is an embedding with only Hopf links (all links have linking number $\pm 1$), shown in Figure~\ref{F:K8hopf}, but this embedding contains 330 links. So it seems there may be a tradeoff between the number of links and the complexity of  the individual links.

\small

\normalsize

\end{document}